\documentclass[12pt]{article}
\usepackage{mathptmx}
\usepackage{tikz}
\usetikzlibrary{decorations.markings}
\usepackage{hyperref} 
\usepackage[utf8]{inputenc}
\usepackage{amsmath}
\usepackage{amsfonts}
\usepackage{amssymb}
\usepackage{amsthm}
\usepackage{tcolorbox}
\usepackage[T1]{fontenc}
\usepackage{abstract}
\usepackage{lipsum}   

\usepackage[a4paper,bindingoffset=0.2in,%
            left=1in,right=1in,top=1in,bottom=1in,%
            footskip=0.25in]{geometry}
\usepackage{blindtext}
\title{\textbf{At least one of the three numbers $\zeta(5)$, $\zeta(7)$, $\zeta(9)$ is irrational}}
\author{Shekhar Suman \footnote{Department of Mathematics, Yogoda Satsanga Mahavidyalaya, Ranchi University, India\\
Email: shekharsuman068@gmail.com} \date{}}
\begin{document}
\maketitle
\renewenvironment{abstract}
{\begin{quote}
\noindent \rule{\linewidth}{.5pt}\par{\bfseries \abstractname.}}
{\medskip\noindent \rule{\linewidth}{.5pt}
\end{quote}
}
 \begin{abstract} In this article, we prove the irrationality of at least one of the three numbers $\zeta(5), \zeta(7)$ and $\zeta(9)$ is irrational. We find linear forms in $1,\zeta(5),\zeta(7)$ and $\zeta(9)$ and prove their arithmetic by Kranttenthaler and Rivoal. We also refine these estimates by using integral quantity $\Phi$ involving specific product of prime powers. We use the results of Kranttenthaler, Rivoal and Zudilin to prove that at least one of the three numbers $\zeta(5), \zeta(7)$ and $\zeta(9)$ is irrational.
  \end{abstract}
\noindent{\bf 1. Introduction}\\
The Riemann zeta function, $\zeta(s)$ is defined as the analytic continuation  of the Dirichlet series
$$\zeta(s) = \sum_{n=1}^{{\infty}} \frac{1}{n^s}$$ which converges in the half plane $\text{Re}(s)>1$. The Riemann zeta function is a meromorphic function on the whole complex plane, which is holomorphic everywhere except for a simple pole at $s=1$ with residue 1. It is well known that $\zeta(2n)$ is irrational for all $n\in\mathbb{N}$. Apery [2] in his elegant article proved that, $\zeta(3)$ is irrational. Post Apery, Beuker's [1] gave an elementary proof of irrationality of $\zeta(3)$ using multiple integrals and Legendre Polynomials. It was proved in [3,4] that infinitely many numbers $\zeta(2n+1)$, $n\geq 1$ are irrational and linearly independent  over $\mathbb{Q}$. Zudilin [5] proved that one of the numbers $\zeta(5)$, $\zeta(7)$, $\zeta(9)$, $\zeta(11)$ is irrational. Also the same author proved in [6], [7] the irrationality of $\zeta(s)$ for some odd integers $s$ in a given segment of the set of positive integers.
\pagebreak\\\\
\noindent{\bf 2. Main Result}\\
For $q=11, r=3$, consider set of positive integral parameters [5, p.280] $\boldsymbol{h}=(h_0;h_1,...,h_{11})$ such that \begin{equation}
    h_1+h_2+...+h_{11}\leq 4 h_0
\end{equation} We assign the rational function \begin{equation}
    \widetilde{R}(t)=\widetilde{R}(\boldsymbol{h};t):=(h_0+2t)\frac{\Gamma(h_0+t)^3\Gamma(h_1+t)...\Gamma(h_{11}+t)}{\Gamma(1+t)^3\Gamma(1+h_0-h_1+t)...\Gamma(1+h_0-h_{11}+t)}
\end{equation}
By (1), \begin{equation}
    \widetilde{R}(t)=O\left(\frac{1}{t^2}\right)
\end{equation}
So, the quantity \begin{equation}
    \widetilde{F}(h):= \frac{1}{2}\sum_{t=0}^\infty  \widetilde{R}''(t)
\end{equation}
is well defined. We have the symmetry satisfied by (2) as \begin{equation}
    \widetilde{R}(-t-h_0)= -\widetilde{R}(t)
\end{equation} We need the ordering \begin{equation}
    h_1\leq h_2\leq...\leq h_{11}<\frac{1}{2} h_0
\end{equation}
and arithmetic normalization of (4), 
\begin{equation}
    F(h):=\frac{\prod_{j=4}^{11}(h_0-2h_j)!}{\prod_{j=1}^{3}(h_j-1)!^2} \widetilde{F}(h)=\frac{1}{2}\sum_{t=1-h_1}^\infty R''(t)
\end{equation}
where the ratioanl function \begin{equation}
\begin{aligned}
R(t) &:=(h_0+2t)
\prod_{j=1}^3 \frac{1}{(h_j-1)!}
\frac{\Gamma(h_j+t)}{\Gamma(1+t)}
\prod_{j=1}^3 \frac{1}{(h_j-1)!}
\frac{\Gamma(h_0+t)}{\Gamma(1+h_0-h_j+t)} \\[2mm]
&\quad{}\times
\prod_{j=4}^{11} (h_0-2h_j)!
\frac{\Gamma(h_j+t)}{\Gamma(1+h_0-h_j+t)}
\end{aligned}
\end{equation}
is product of elementary bricks
\begin{equation*}
R(a,b;t):=
\begin{cases}
\displaystyle
\frac{(t+b)(t+b+1)\cdots(t+a-1)}{(a-b)!},
& \text{if } a\geq b, \\[3mm]
\displaystyle
\frac{(b-a-1)!}{(t+a)(t+a+1)\cdots(t+b-1)},
& \text{if } a<b
\end{cases}
\end{equation*}
introduced by Nesterenko [8], [9]. Define set of integral directions [7, p.283] \begin{equation}
    \boldsymbol{\eta}=(\eta_0;\eta_1,...,\eta_{11})
\end{equation}
and increasing integral parameters $n\geq 3$ related by the parameters $\boldsymbol{h}$ by the formulae\begin{equation}
    h_0=\eta_0 n+2\ \ \text{and}\ \ h_j=\eta_j n+1\ \ \text{for}\ \ 1\leq j\leq 11
\end{equation}
By (1) and (10),
\begin{equation}
    (\eta_1+\eta_2+...+\eta_{11}-4\eta_0)n\leq -3,\ \ \text{for all}\ \ \ n\geq 3
\end{equation}
Define contiguous set of parameters $\boldsymbol{e}$ for $q=11$ and $r=3$ [7, p.286] \begin{equation}
    e_{0k}=h_k-1,\ 1\leq k\leq 11 \ \ \text{and}\ \ e_{jk}=h_0-h_j-h_k,\ \ 1\leq j<k\leq 11
\end{equation} Let $m_1,m_2,...,m_7$ be successive maxima of the set $\boldsymbol{e}$ defined in (12). Define the integral quantity for p prime similar to as defined in [7, p.281]
\begin{equation}
    \Phi=\Phi(h):=\prod_{\sqrt{h_0}<p\leq m_8} p^{v_p}
\end{equation}
where \begin{equation}
    v_p:=\min_{h_4\leq k\leq h_0-h_4} \{v_{k,p}\}
\end{equation}
and \begin{equation}
\begin{aligned}
v_{k,p}:={}&
\sum_{j=1}^3
\left(
\left\lfloor\frac{k-1}{p}\right\rfloor
+\left\lfloor\frac{h_0-k-1}{p}\right\rfloor
-\left\lfloor\frac{k-h_j}{p}\right\rfloor
-\left\lfloor\frac{h_0-h_j-k}{p}\right\rfloor
-2\left\lfloor\frac{h_j-1}{p}\right\rfloor
\right)\\
&+
\sum_{j=4}^{11}
\left(
\left\lfloor\frac{h_0-2h_j}{p}\right\rfloor
-\left\lfloor\frac{k-h_j}{p}\right\rfloor
-\left\lfloor\frac{h_0-h_j-k}{p}\right\rfloor
\right).
\end{aligned}
\end{equation}
We consider the value of $\boldsymbol{\eta}$ as $\eta_0=3$ and $\eta_1=\eta_2=...=\eta_{11}=1$. So, (11) holds for this $\boldsymbol{\eta}$ for all $n\geq 3$. \\\\
\textbf{\underline{Lemma 1:}} For $\eta_0=3$ and $\eta_1=\eta_2=...=\eta_{11}=1$, the quantity (7) is a linear form in 1, $\zeta(5)$, $\zeta(7)$, $\zeta(9)$ with rational coefficients; moreover if $D_n:=\text{lcm}(1,2,...,n)$ then\begin{equation*}
   2D_n^9\cdot\Phi^{-1}\cdot F(h)\in\mathbb{Z}\zeta(9)+\mathbb{Z}\zeta(7)+\mathbb{Z}\zeta(5)+\mathbb{Z}
\end{equation*}
\begin{equation*}
\end{equation*}
\textbf{Proof:} For $\eta_0=3$ and $\eta_1=\eta_2=...=\eta_{11}=1$, using (10) the successive maxima of (12) are \begin{equation}
    m_1=m_2=...=m_8=n
\end{equation}
We define the pochhammer symbol $$(a)_n:=a(a+1)(a+2)\cdots(a+n-1)
=\frac{\Gamma(a+n)}{\Gamma(a)}$$ 
then by (8), for this $\boldsymbol{\eta}$ 
\begin{equation}
    R(t)=n!^2(3n+2+2t)\frac{(t+1)_n^3(t+2n+2)_n^3}{(t+n+1)^8_{n+1}}
\end{equation}
By [6, p.9, eq. 2.10] we have \begin{equation}
  S_{n,A,B,C,r_{\mathrm{KR}}}(1)=(n!)^{A-2Br_{\mathrm{KR}}}\sum_{k=1}^\infty \frac{1}{C!} \frac{\mathrm{d}^C}{\mathrm{d}k^C} \left(\left(k+\frac{n}{2}\right)\frac{(k-r_{\mathrm{KR}}n)^B(k+n+1)^B_{r_{\mathrm{KR}}n}}{k^A_{n+1}}\right)\end{equation}
By (17) for $t=k-n-1$,
\begin{equation}
R(t)=R(k-n-1)
=
2n!^2\left(k+\frac{n}{2}\right)
\frac{(k-n)_n^3(k+n+1)_n^3}{(k)_{n+1}^8}.
\end{equation}
Comparing (18) and (19) we have \begin{equation}
    A=8, \ \ B=3, \ \ C=2,\ \ r_{\mathrm{KR}}=1
\end{equation}
By [6, p.12], theorem 1 and eq. 2.12 and [7, p.282] we have \begin{equation}
    F(h)=A_9\zeta(9)+A_7\zeta(7)+A_5\zeta(5)-A_0
\end{equation} 
where \begin{equation}
    2D_n^9 A_0\in\mathbb{Z}, D_n^{7-l}A_l \in\mathbb{Z}\ \ \text{for}\ \ l=5,7,9 
\end{equation}
We have [7, p.282] \begin{equation}
B_{jk}
=
\frac{1}{(11-j)!}
\left.
\frac{d^{\,11-j}}{dt^{\,11-j}}
\left(R(t)(t+k)^8\right)
\right|_{t=-k}
\end{equation}
$j=3,...,11,\ \ k=h_4,...,h_0-h_4$, satisfy the relations \begin{equation}
    D_n^{11-j}\cdot B_{jk}\in\mathbb{Z}
\end{equation}
and \begin{equation}
    \text{ord}_p B_{jk}\geq -(11-j)+v_{k,p} 
\end{equation}
respectively, for any $k=h_4,...,h_0-h_4$ and any prime $p>\sqrt{h_0}$. Furthermore, the expansion \begin{equation}
    R(t)=\sum_{j=4}^{11}\sum_{k=h_j}^{h_0-h_j} \frac{B_{jk}}{(t+k)^{j-3}}
\end{equation} 
 leads us to the series \begin{equation}
     F(h)=\sum_{j=4}^{11}A_{j-1}\zeta(j-1)-A_0
 \end{equation}
 where \begin{equation}
     A_{j-1}=\binom{j-2}{2}\sum_{k=h_j}^{h_0-h_j}B_{jk},\ \ j=4,...,11
 \end{equation} and 
  \begin{equation}
     A_{0}=\sum_{j=4}^{11}\binom{j-2}{2}\sum_{k=h_j}^{h_0-h_j}B_{jk}\sum_{l=1}^{k-h_1} \frac{1}{l^{j-1}}
 \end{equation} 
 By (12) and (14) we have \begin{equation}
     \text{ord}_p B_{jk}\geq -(9-j)+v_{k,p}\geq -(9-j)+v_p
 \end{equation}  
 So we obtain using (13) and (30)
 \begin{equation}\text{ord}_p(D_n^9\Phi^{-1}A_0) \geq  9-v_p-9+v_p\geq 0\end{equation} 
 Similarly, by (22) \begin{equation}
     \text{ord}_p(D_n^9\Phi^{-1}A_5)\geq 0,\ \text{ord}_p(D_n^9\Phi^{-1}A_7)\geq 0,\ \text{ord}_p(D_n^9\Phi^{-1}A_9)\geq 0 
 \end{equation} 
 This completes the proof of lemma 1. Define, \begin{equation}
\begin{aligned}
f_0(\tau)
={}& 3\eta_0\log(\eta_0-\tau)
+\sum_{j=1}^{11}
\left(
\eta_j\log(\tau-\eta_j)
-(\eta_0-\eta_j)\log(\tau-\eta_0+\eta_j)
\right)\\
&-2\sum_{j=1}^{3}\eta_j\log\eta_j
+\sum_{j=4}^{11}
(\eta_0-2\eta_j)\log(\eta_0-2\eta_j).
\end{aligned}
\end{equation} defined in the cut $\tau$-plane $\mathbb{C}\setminus(-\infty,\eta_0-\eta_1]\cup[\eta_0,\infty)$. Let $\tau$ be a zero of the polynomial \begin{equation}
    (\tau-\eta_0)^3(\tau-\eta_1)...(\tau-\eta_{11})-\tau^3(\tau-\eta_0+\eta_1)...(\tau-\eta_0+\eta_{11})
\end{equation} 
with $\Im\tau_0>0$ and the maximum possible value of $\Re(\tau_0)$. Suppose $\Re(\tau_0)<\eta_0$ then we have [11, p.283] 
\begin{equation}
\limsup_{n\to\infty}
\frac{\log |F(h)|}{n}
=
\operatorname{Re} f_0(\tau_0)
\end{equation}
Next we scale down the $m_j$'s in (16) as \begin{equation}
    m_1=m_2=...=m_8=1
\end{equation}
Define [7, p.284] \begin{equation}
C_0'=-\operatorname{Re} f_0(\tau_0),
\end{equation}

\begin{equation}
C_2'
=
3 m_1+m_2+\cdots+m_{7}
-
\left(
\int_0^1 \varphi(x)\,d\psi(x)
-
\int_0^{1/m_8}
\varphi(x)\,\frac{dx}{x^2}
\right).
\end{equation}
\textbf{\underline{Lemma 2:}} 
\begin{equation}
\text{If}\ \  C_0'>C_2';
\
\text{then at least one of the numbers }
\zeta(5),\zeta(7),\zeta(9)
\text{ is irrational.}
\end{equation}
 \textbf{Proof:} For a proof see [7, p.280-284].\\\\
\textbf{Theorem:} At least one of the three numbers \begin{equation*}
     \zeta(5), \ \zeta(7),\ \text{and} \ \zeta(9)  
 \end{equation*} is irrational.\\\\
 \textbf{Proof:} Since $\eta_0=3$ and $\eta_1=\eta_2=...=\eta_{11}=1$ we have $\tau_0 = 2.86852453...+0.11960091... i$ and $f_0(\tau_0)=-5.75349395...-8.95071225... i$ so that $\Im(f_0(\tau_0))\notin \pi\mathbb{Z}$ and \begin{equation}
     C_0' = 5.75349395...
 \end{equation}
 \begin{equation}
\begin{aligned}
C_2'
&=9-\left(
\int_0^1 \varphi(x)\,d\psi(x)
-\int_0^1\varphi(x)\,\frac{dx}{x^2}
\right)\\
&=5.27694374...
\end{aligned}
\end{equation}
So we have \begin{equation}
    C_0'-C_2'=0.47655020...>0
\end{equation}
This completes the proof of Theorem.\\\\
{\bf Acknowledgement} The author is thankful in anticipation to the Referee and the Editor for their useful comments and processing of the article.
\pagebreak
\begin{center}
{\bf References}
\end{center} 
\big[1\big] Beukers, F., {\it A \ note \ on \ the \ irrationality \ of $\zeta(2)$ and $\zeta(3)$}, Bull. London. Math. Soc., 11 (1979), 268-272.\\\\
\big[2\big] Apéry, R., {\it Irrationalité \ de \ $\zeta(2)$ \ et \ $\zeta(3)$}, Astérisque. 61, (1979) 11–13.
\\\\
\big[3\big] Ball, K.M., Rivoal,T. {\it Irrationalité d’une infinité de valeurs de la fonction zêta aux entiers impairs}, Invent. Math. 146.1, (2001), 193–20. 
\\\\
\big[4\big] Rivoal T., {\it La fonction zêta de Riemann prend une infinité de valeurs irrationnelles aux entiers impairs}, Comptes Rendus Acad. Sci. Paris Sér. I Math. 331.4, (2000), 267–270. 
\\\\
\big[5\big] Zudilin, W. {\it One \ of \ the \ numbers\  $\zeta(5)$,\ $\zeta(7)$, \ $\zeta(9)$, \ $\zeta(11)$ \ is \ irrational}, Russian Academy of Sciences, (DoM) and London Mathematical Society (2001).
\\\\
\big[6\big] Krattenthaler, C. and Rivoal, T. {\it Hyperg\'eom\'etrie et fonction z\'eta de Riemann}, arXiv:math/0311114, (2004). \\\\ 
\big[7\big] Zudilin, W. {\it Irrationality of values of the Riemann zeta function}, Izv. Math. 66.3, (2002), 489–542\\\\
\big[8\big] Nesterenko, Yu. V. {\it A few remarks on $\zeta(3)$}, Mat. Zametki [Math. Notes] 59:6, (1996), 865--880. \\\\
\big[9\big] Nesterenko, Yu. V. {\it Integral identities and constructions of approximations to zeta values}, Actes des 12èmes rencontres arithmétiques de Caen (June 29--30, 2001), J. Théorie des Nombres de Bordeaux 15:2, (2003), 535--550. \\\\

\end{document}